\documentclass{article}
\usepackage{authblk}  
\usepackage{marvosym} 
\usepackage{graphicx} 
\usepackage{xcolor}   

\graphicspath{{\subfix{../images/}}}
\usepackage{blindtext}
\usepackage{caption}
\usepackage{subcaption}
\usepackage{amsthm}
\usepackage{hyperref}
\usepackage[margin=1in]{geometry}

\theoremstyle{definition}
\newtheorem{definition}{Definition}

\newcommand{\V}{{\cal V}}
\newcommand{\E}{{\cal E}}
\newcommand{\sub}{\subseteq}

\title{Dynamical Properties of Random Boolean Hypernetworks}

\author[1(\Letter)]{Kevin M. Stoltz}
\author[1, 2]{Cliff A. Joslyn}
\affil[1]{Systems Science and Industrial Engineering Department, Binghamton University (State University of New York), Binghamton, USA}
\affil[ ]{\texttt{\{kstoltz1, cajoslyn\}@binghamton.edu}}
\affil[2]{Pacific Northwest National Laboratory, Seattle, USA}
\affil[ ]{\texttt{cliff.joslyn@pnnl.gov}}
\date{\today}

\begin{document}
\maketitle

\begin{abstract}
Boolean networks are a valuable class of discrete dynamical systems models, but they remain fundamentally limited by their inability to capture multi-way interactions in their components. To remedy this limitation, we propose a model of Boolean hypernetworks, which generalize standard Boolean networks. Utilizing the bijection between hypernetworks and bipartite networks, we show how Boolean hypernetworks generalize standard Boolean networks. We derive ensembles of Boolean hypernetworks from standard random Boolean networks and simulate the dynamics of each. Our results indicate that several properties of Boolean network dynamics are affected by the addition of multi-way interactions, and that these additions can have stabilizing or destabilizing effects.   
\end{abstract}

\section{Introduction} \label{intro}
Random Boolean networks, introduced by Kauffman \cite{kauffman_metabolic_1969} in an effort to model genetic interactions, can be characterized by three key attributes: (1) directed networks constructed with Boolean state variables on vertices, random edge placement with fixed in-degree $k$ for each vertex, and random state update rules produced a surprising amount of order when their dynamics were simulated. After a sufficient amount of time, they entered stable cycles of states termed attractors that comprised a small fraction of the total state space; (2) attractor cycles went from ordered, characterized by small cycles and few attractors, to chaotic, characterized by much longer cycles, as $k$ was increased; and (3) networks with $k = 2$ were found to represent the critical point between ordered ($k < 2$) and chaotic ($k > 2$) dynamics.

Boolean network models are still widely used, often to model biological networks in which connections and state update rules are the result of literature reviews and experimental evidence instead of random assignment \cite{manicka_effective_2022, gates_effective_2021, gomez_tejeda_zanudo_network_2017}. While such models have yielded valuable insight for basic and applied research, they suffer from a fundamental limitation inherent to all network-based methods: they are limited to modeling pairwise interactions between vertices. 

Hypernetworks are a class of models that generalize networks by allowing for an arbitrary number of vertices, bounded above by the total number of vertices in the hypernetwork, to be included in each edge. For this reason, edges in hypernetworks are called hyperedges, and they capture the higher order interactions missing in network science but ubiquitous in real-world complex systems \cite{braha_hypernetwork_2021, klamt_hypergraphs_2009}. Evidence for the importance of taking these higher order interactions into account can be seen from the fact that centrality measures on hypernetworks have been shown to outperform standard graph centrality measures in identifying genes that play a role in pathogenic responses to viral infections \cite{feng_hypergraph_2021}. As hypernetworks are generalizations of networks, every network is in fact a hypernetwork in which each edge contains exactly two vertices, i.e. a two-uniform hypernetwork. Furthermore, the directed networks that form the structure of Kauffman's random Boolean networks find a generalization in directed hypernetworks \cite{gallo_directed_1993, brenner_exchange_2017}. Crucially, hypernetworks are bijective to bipartite graphs in which the two disjoint sets of bipartite network vertices are the hyperedges and (hyper)vertices of the hypernetwork.

Our aim in this work is to generalize random Boolean networks by defining and characterizing a model of random Boolean hypernetworks. We believe that the study of such models could yield valuable insights into the dynamics of higher order networks. The remainder of this work is organized as follows: after giving an overview of related work in section \ref{related} we formally introduce random Boolean networks and their bipartite representations in section \ref{prelim}. In section \ref{rbh}, we provide a small example of the transformation of a random Boolean network to its isomorphic bipartite network, and then its generalization to a random Boolean hypernetwork. In section \ref{results}, we simulate ensembles of random Boolean hypernetworks and report our results. We conclude with a discussion in section \ref{discussion}.

Throughout this work, we make use of the bijective relationship between bipartite graphs and hypernetworks. All hypernetworks, whether they be two-uniform or composed of multi-way interactions, are constructed and simulated in their bipartite form using the NetworkX package \cite{hagberg_exploring_2008}. When analyzing attractors, we first generate the state transition graph from simulation results, then use the PyBoolNet package to calculate attractor period lengths \cite{klarner_pyboolnet_2017}.

\section{Related Work} \label{related}
While we have yet to find explicit descriptions of Boolean hypernetworks in the literature, there does exist an active and growing body of research on the topic of bipartite Boolean networks. Introduced by Graudenzi in an effort to faithfully represent gene-protein interactions, this model was comprised of G nodes (representing genes) and P nodes (representing proteins) \cite{graudenzi_new_2009}. It is characterized by a decay phase parameter on P nodes to account for protein degradation in the cytoplasm, and the topological constraint that G node out-degree must equal 1 since each genes encodes a single protein. Dynamical properties of this model have been examined with respect to changes in the decay phase parameter on P nodes \cite{graudenzi_dynamical_2011, graudenzi_robustness_2011, sapienza_dynamical_2018}, but the out-degree constraint on G nodes has crucial implications. The constraint on G node out-degree leads to G and P nodes always being present in equal proportions. Furthermore, each G-node having out-degree 1 implies that each P node has in-degree 1. The models presented in our work loosen both of these constraints: the two subsets of bipartite vertices that comprise our networks can have different cardinalities; and while one subset of our bipartite vertices is constrained by their out-degree being equal to 1, their successors are \textit{not} constrained by having in-degree equal to 1. 

A more recent model proposed in \cite{hannam_percolation_2019} and \cite{torrisi_percolation_2020} used transcription factors instead of single genes as P nodes. As transcription factors can be multi-protein complexes, this model looses both the G node and P node degree constraints described above. Crucially, the activation status of the P nodes modeling transcription factors are solely determined by the logical AND function since multi-protein complexes require each subunit to be expressed. Similar constraints on logical functions are used in \cite{lemke_essentiality_2004, ghim_lethality_2006, lee_branching_2012}, which use logical AND and OR exclusively to model metabolic networks in \textit{Escherichia coli}. Our model loosens these constraints by assigning to each node a randomly-chosen Boolean function of the states of its inputs.

The introduction of gene-transcription factor networks opened up a novel set of questions about the number of different ways one gene can influence another when taking transcription factor intermediaries into account. An answer to this question can be found by studying the composition of Boolean functions, and \cite{fink_boolean_nodate} gave a method for determining the number of distinct Boolean functions given a particular compositional structure. \cite{yadav_relative_2022} further characterized these functional composition structures, and examined the enrichment of different compositional structures in real biological networks.

\section{Preliminaries} \label{prelim}

\subsection{Hypergraphs and Bipartite Graphs}

Hypergraphs (directed and undirected) are isomorphic to corresponding bipartite graphs (directed and undirected). 

\begin{definition}[Undirected Hypergraph]
    A hypergraph $H = ({\cal V}, {\cal E})$ is a set of vertices ${\cal V} = \{v_i\}, 1 \le i \le n$, and an indexed family of hyperedges ${\cal E} = \{e_j\}, 1 \le j \le m$, such that $e_j \subseteq {\cal V}$ for all $e_j \in {\cal E}$.
\end{definition}

We recognize graphs as a special case of hypergraphs in which each hyperedge has cardinality two. Furthermore, every hypergraph $H=({\cal V},{\cal E})$ maps  bijectively to a bipartite graph ${\cal G} = ({\cal V} \sqcup {\cal E}, {\cal D})$ with new vertex set ${\cal V} \sqcup {\cal E}$ and edges ${\cal D} \subseteq {{\cal V} \cup {\cal E} \choose 2}$ such that $d=\{v_i,e_j\} \in {\cal D}$ iff $\exists e_j \in {\cal E}, v_i \in e_j$. Note that thereby ${\cal G}$ is bipartite with two parts ${\cal V}$ and ${\cal E}$ such that for all edges $d = \{x_k,x_{k'}\} \in {\cal D}$, for $x_k,x_{k'} \in {\cal V} \sqcup {\cal E}, 1 \le k,k' \le n+m$, neither $d \subseteq {\cal V}$ nor $d \subseteq {\cal E}$. See \cite{aksoy_hypernetwork_2020} for more details.

\begin{figure}[ht]
\centering
\includegraphics[scale=0.23]{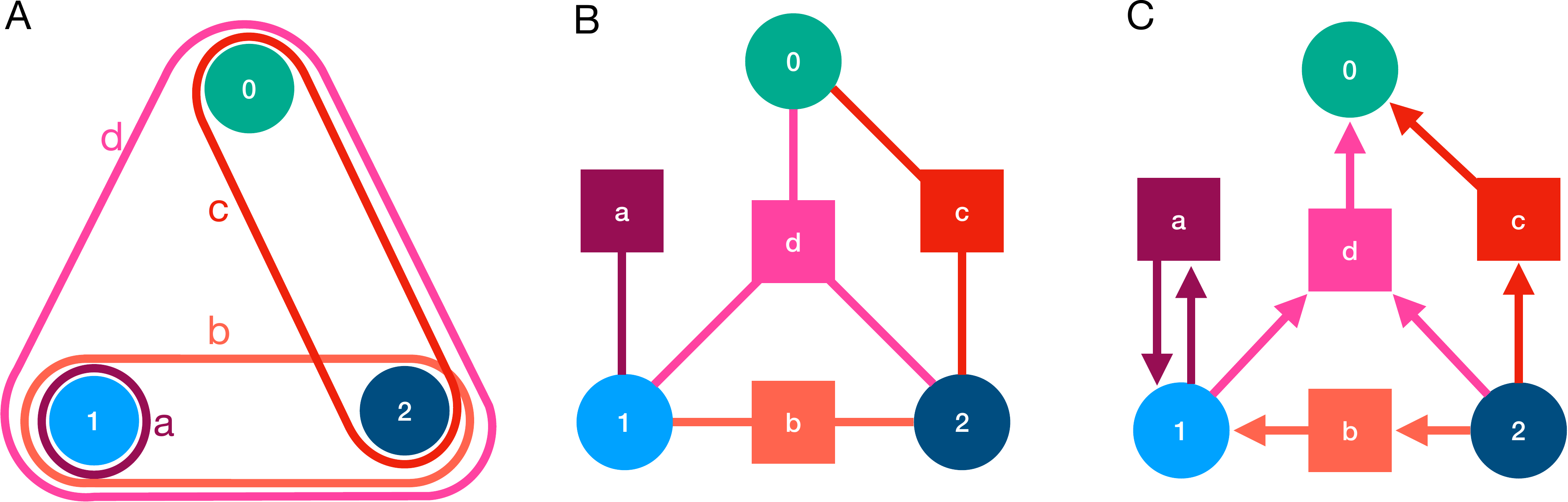}
\caption{(A) and (B) A small example of an undirected hypergraph and its bipartite representation, respectively. (C) The bipartite representation of one possible directed hypergraph that can be derived from the same vertex set as in B.}
\label{fig:hypergraph relationships}
\end{figure}

A simple example of an undirected hypergraph is shown in Figure \ref{fig:hypergraph relationships}A, as an ``Euler diagram'', with hyperedges $e_j$ surrounding their vertices $v_i$.  Figure \ref{fig:hypergraph relationships}B shows the corresponding bipartite graph, now with the circles as the $v_i$ part and squares as the $e_j$ part. The edges $d \in {\cal D}$ are then shown as undirected graph edges connecting the two parts.

Hypergraphs also have a directed form, analogously generalizing directed graphs.
\begin{definition}[Directed Hypergraph]
A hypergraph $H=(\V,\E)$ is directed when each edge $e_j \sub \V$ is represented as a cover consisting of a head $e^h_j \sub e_j$ and a tail $e^t_j \sub e_j$ such that $e^h_j \cup e^t_j = e_j$.
\end{definition}

Similar to the case of undericted hypergraphs, directed hypergraphs map bijectively to directed bipartite graphs. The directed hypergraphs that underlie our proposed Boolean hypernetworks differ from those introduced by \cite{gallo_directed_1993} in that we do not require the vertex subsets that compose a single directed hyperedge to be disjoint, as indicated above, to allow for self-loops in our construction. In other words, it may be that $e^h_j \cap e^t_j \neq \emptyset$, with vertices in the intersection represented as a 2-cycle in the corresponding directed bipartite graph. 

Figure \ref{fig:hypergraph relationships}C gives the bipartite representation of one possible directed hypergraph that can be derived from \ref{fig:hypergraph relationships}B; that is, one possible division of the hyperedges into heads and tails. Note that we are not showing an Euler diagram equivalent to \ref{fig:hypergraph relationships}A in the undirected case, since these are extremely challenging to visualize in that form, despite their formal equivalence.

We will also use the following notation below for directed bipartite graphs.

\begin{definition}[Directed Bipartite Graph]
    A directed bipartite graph $\mathcal{G} = (\mathcal{V}, \mathcal{E}, \mathcal{D})$ is composed of: 
    \begin{enumerate}
        \item Vertex subsets $\mathcal{V} = \{v_i\}$, $1 \leq i \leq n$, each with in-degree $k$ and $\mathcal{E} = \{e_j\}$, $n + 1 \leq j \leq n + m$, each with in-degree = out-degree = 1 such that $\mathcal{V} \cap \mathcal{E} = \emptyset$; and
    
        \item a set of directed edges $\mathcal{D} \subseteq (\mathcal{V} \sqcup \mathcal{E}) \times (\mathcal{V} \sqcup \mathcal{E})$ such that $\forall v_p, v_q \in \mathcal{V}$ and $e_{p'}, e_{q'} \in \mathcal{E}$, we have that $(v_p, v_q), (v_q, v_p), (e_{p'}, e_{q'}), (e_{q'}, e_{p'}) \notin \mathcal{D}$. 
    \end{enumerate}
\end{definition}

Note that now by construction, $|\mathcal{E}| = m$.   

\subsection{Random Boolean Networks}

Random Boolean networks were first introduced by Kauffman as a way to model biological genetic networks \cite{kauffman_metabolic_1969}. Genes and their interactions are modeled as directed graphs where each vertex takes on a binary state that can be influenced by its incoming connections. While other network topologies, including scale-free topologies \cite{aldana_boolean_2003} have since been studied, Kauffman's original model focused on the case where each vertex has fixed in-degree $k$. Once the size of the network and $k$ have been specified, connections are placed randomly and each vertex is randomly assigned one of $2^{2^k}$ possible Boolean state update functions. 

\begin{definition}[Boolean Network] \label{def:BN}
A random Boolean network $BN = (G, S, F)$ is composed of:
\begin{enumerate}
    
\item A directed graph $G = (V, E)$ composed of a finite, non-empty set of $n$ vertices $V = \{v_i\}, 1 \leq i \leq n$, each with in-degree $k$, and a set of $m = nk$ directed edges $E = \{e_j\}, 1 \leq j \leq m$ where $e_j \in V \times V$ for all $e_j$;

\item A set of state variables $S = \{s_i\}$; and

\item State transition functions $F = \{f_i\}, 1 \leq i \leq n$, each associated with $V$. 
 
\end{enumerate}

\end{definition}
Specifically, each vertex $v_i \in V$ has a state variable $s_i \in \{0, 1\}$ whose time evolution is governed by state transition function $f_i: \{0, 1\}^k \rightarrow \{0, 1\}$. Letting $s_i^t$ denote $s_i$ at time $t$, the temporal evolution of $s_i$ is described by $s_i^{t+1} = f_i(s_{w1}^t, ..., s_{wk}^t)$ where $(v_{w1}, v_i), ..., (v_{wk}, v_i) \in E$. The vertices $v_{w1}, ..., v_{wk}$ are the inputs to vertex $v_i$. Finally, let $S^t = \{s^t_1, ..., s^t_n\}$ be the set of states of all vertices at time t. 

A simple example of a $BN$ for $n=3$ and $k=2$ is shown in Figure \ref{fig:network relationships}A. To generate the bipartite graph isomorphic to this network, we recall that each directed graph is a 2-uniform directed hypergraph where each hyperedge contains one in vertex and one out vertex. Furthermore, each directed hypergraph is isomorphic to a directed bipartite graph. It follows that every $BN$ is isomorphic to a directed bipartite network. 

\begin{definition}[Bipartite Random Boolean Network]
    A bipartite random Boolean network $BBN = (\mathcal{G}, \mathcal{S}, \mathcal{F})$ is composed of 
    \begin{enumerate}
        \item a directed bipartite graph $\mathcal{G}$;

        \item sets of state variables $\mathcal{S} = \{\mathcal{S}_\mathcal{V} = \{s_i \in \{0, 1\}\}, \mathcal{S}_\mathcal{E} = \{s_j \in \{0, 1\}\}\}$ and state transition functions $\mathcal{F} = \{\mathcal{F}_\mathcal{V} = \{f_i: \{0, 1\}^k \rightarrow \{0, 1\}\}, \mathcal{F}_\mathcal{E}\ = \{f_j: \{0, 1\} \rightarrow \{0, 1\}\}\}$, $1 \leq i \leq n$, $n+1 \leq j \leq n+m$ associated with $\mathcal{V}$ and $\mathcal{E}$, respectively.
    \end{enumerate}
\end{definition}

Then we can denote $s^t_x$ as $s_x$ at time $t$ and describe the temporal evolution of $s_x$ as in Definition \ref{def:BN}. To define the bipartite network isomorphic to a random Boolean hypernetwork, we modify two properties of a $BBN$.

\begin{definition} [Bipartite Random Boolean Hypernetwork]
    A bipartite random Boolean hypernetwork $BH$ extension of a $BBN$ is identical to the $BBN$ except in two properties: 
    \begin{enumerate}
        \item the in-degree of each vertex in $\mathcal{E}$ is $l \ge 1$;

        \item $\mathcal{F}_\mathcal{E}\ = \{f_j: \{0, 1\}^l \rightarrow \{0, 1\}\}$.
    \end{enumerate}
\end{definition}

In the case where $l = 1$, $BH = BBN$.

\section{Random Boolean Hypernetworks} \label{rbh}
\begin{figure}[ht]
\centering
\includegraphics[scale=0.23]{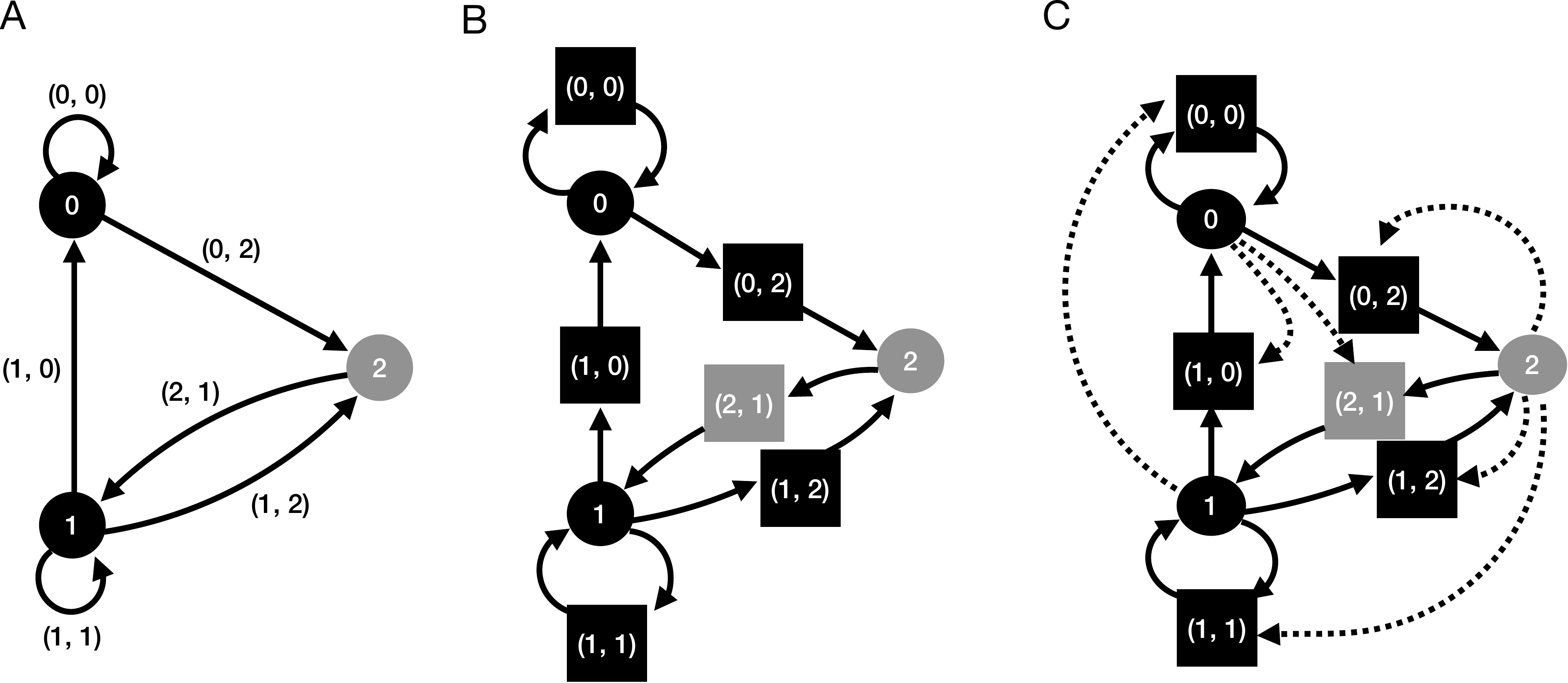}
\caption{Deriving a random Boolean hypernetwork from a random Boolean network. (A) A $BN$ with $k=2$. Vertex colors denote state variables with black indicating 1 and gray indicating 0. (B) The $BBN$ derived from A. The vertex (2, 1) takes the state of its single input, as specified by the identity update rule. (C) A $BH$ extending B, where dashed edges are additional edges required when increasing $l$ from 1 to 2. Here, the vertex (2, 1) may, or may not take on the same state value as vertex 2.}
\label{fig:network relationships}
\end{figure}

The $BN$ in Figure \ref{fig:network relationships}A is isomorphic to the $BBN$ in Figure \ref{fig:network relationships}B. To generate Figure \ref{fig:network relationships}B from Figure \ref{fig:network relationships}A, we proceed with the following steps: (1) For all $v_i \in V$, assign $v_i$ to $v_i \in \mathcal{V}$; (2) for $e_j \in E$ composed of $(v_p, v_q)$, assign $e_j$ to $e_j \in \mathcal{E}$; (3) construct two edges: $(v_p, e_j)$ and $(e_j, v_q)$; (4) for $f_i \in F$, assign $f_i$ to $f_i \in \mathcal{F}_\mathcal{V}$; (5) for $f_j \in \mathcal{F}_\mathcal{E}$, assign $f_j$ the Boolean identity function, in which the output value equals the state of the input value.

As an example, all vertices in Figure \ref{fig:network relationships}A are isomorphic to a circular vertex in Figure \ref{fig:network relationships}B. Furthermore, the edge $(0, 2)$ in Figure \ref{fig:network relationships}A is isomorphic to vertex $(0, 2)$ in Figure \ref{fig:network relationships}B. Vertex $(0, 2)$ in Figure \ref{fig:network relationships}B has a single incoming edge from vertex $0$ and a single outgoing edge to vertex $2$. It follows from this construction that in-degree = $k$ for all $v_i \in \mathcal{V}$ and in-degree = out-degree = 1 for all $e_j \in \mathcal{E}$. This construction also implies that if $S^0 = S^0_\mathcal{V}$, simulation of the networks in Figure \ref{fig:network relationships}A and \ref{fig:network relationships}B will yield identical state space trajectories.

\begin{figure}[ht]
\centering
\includegraphics[scale=0.225]{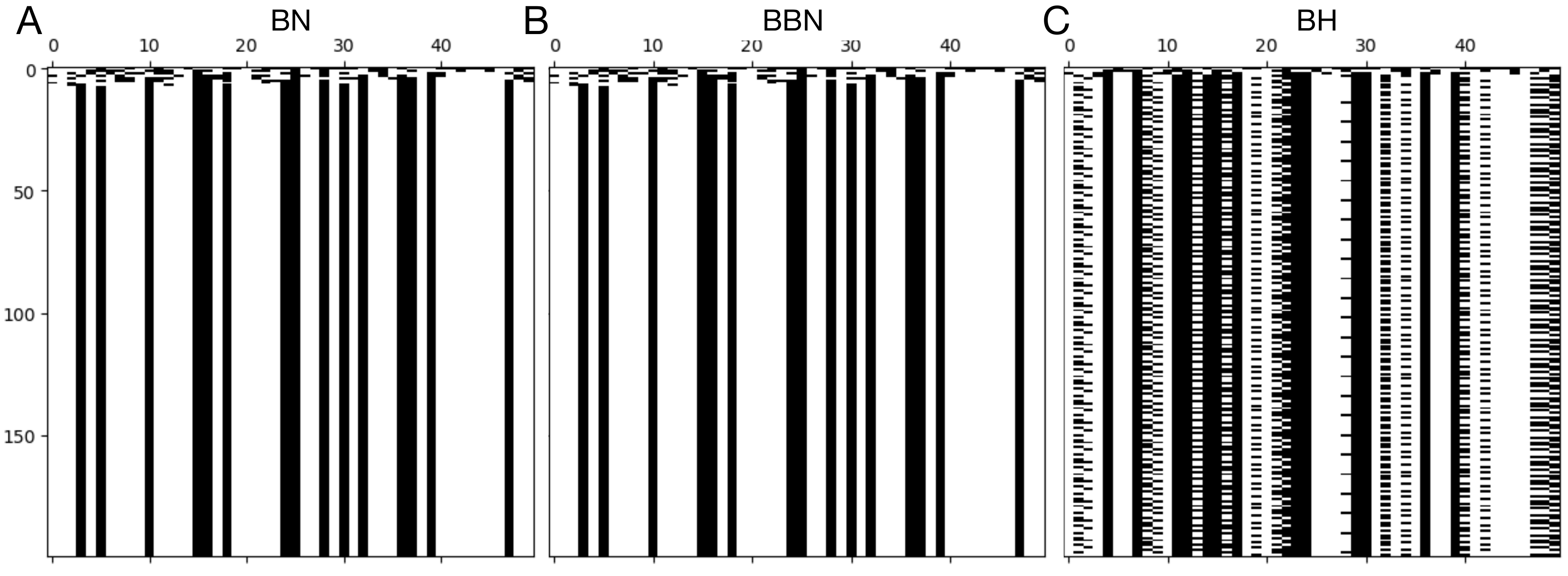}
\caption{State space trajectories for a $BN$ of $n=50$ and the $BBN$ and $BH$ derived from it.  (A) $BN$ state space trajectory. Each square is the state of a vertex (x-axis) at a specified time point (y-axis). Each row is the state of the network at a given time, while each column is the series of state transitions for a single vertex. (B) The exact trajectory of A is recovered when each $\mathcal{V}$ vertex of the $BBN$ is given the same initial condition as its corresponding vertex in the $BN$. (C) The same initial conditions give rise to a state space trajectory with a different attractor period length for the $BH$ derived from the $BBN$ in B.}
\label{fig:transition dynamics}
\end{figure}

$BBNs$ are simulated with an updating scheme in which $\mathcal{V}$ and $\mathcal{E}$ are updated on alternating time steps. This means that for an $BN$ simulated for 100 time steps, its isomorphic $BBN$ would require 200 time steps in order for $\mathcal{V}$ to recieve the same number of updates as $V$. To ensure that a $BN$ and its isomorphic $BBN$ with identical initial states produce identical state space dynamics, we constructed 200 different $BNs$ for $n=50$ and $k=2$ and derived their $BBNs$. After ensuring that initial states were identical, we simulated each $BN (BBN)$ for 100 (200) time steps. In all cases, identical states space trajectories were achieved. An example is shown in Figure \ref{fig:transition dynamics} (A and B). Note that for the $BBN$, only the state space trajectories of $\mathcal{V}$ vertices on alternating time steps is shown, as the $BN$ has no corresponding $\mathcal{E}$ vertices.

$BBNs$ are, by definition, constrained by two facts: (1) all vertices in $\mathcal{E}$ must have in-degree = out-degree = 1; and (2) all vertices in $\mathcal{E}$ must take the identity update rule. In order to construct bipartite networks isomorphic to Boolean hypernetworks ($BH$), we remove these two constraints. In fact, removal of the first constraint necessitates removal of the second as extra inputs to vertices in $\mathcal{E}$ imply that the identity Boolean function can no longer be used. We proceed with the following steps: (1) choose a value $l > 1$ - this value specifies the in-degree for each vertex in $\mathcal{E}$. (2) for each $e_j \in \mathcal{E}$, randomly add inputs to $e_j$ from vertices in $\mathcal{V}$; and (3) randomly assign one of the $2^{2^l}$ possible state update functions. 

An example is shown in Figure \ref{fig:network relationships}C where dashed lines denote newly-added edges. Obviously, the network in Figure \ref{fig:network relationships}C would not be expected to have the same state space dynamics as the $BBN$ in Figure \ref{fig:network relationships}B. Not only do we observe this, we have observed examples in which the dynamics differ to such an extent that attractor length dynamics differ between a $BBN$ and the $BH$ derived from it (compare Figure \ref{fig:transition dynamics}B and Figure \ref{fig:transition dynamics}C). Given these observations, we decided to investigate the extent to which adding inputs to $\mathcal{E}$ vertices alters state space trajectories.

Firstly, we examined the trajectory overlap of vertex states between $BBNs$ and their derived $BHs$. A measure of overlap between two network configurations was introduced in \cite{aldana_boolean_2003}. We make slight modifications to define the state space trajectory overlap between a $BBN$ and a $BH$ as follows. Given a $BBN$ and the $BH$ that extends it, let the state space trajectory overlap between the two networks be

\begin{equation}
    O = \frac{1}{(n+m)T} \sum_{t=1}^{T} \sum_{z=1}^{n+m} |s^{bt}_z - s^{ht}_z|, 
\end{equation}
where $s^{bt}_z$ and $s^{ht}_z$ are states for single $BBN$ and $BH$ vertices at time $t$, respectively. When $O = 1$, the two state space trajectories are identical, whereas $O = 0.5$ implies nearly totally independent trajectories.

In order to determine how changes in $l$ affected trajectory overlap, we constructed 200 $BNs$ with $n=50$, then derived their corresponding $BBNs$ and $BHs$. Because increases in $k$ have already been shown to move $BNs$ from stable to chaotic dynamics, values of $k$ remained fixed as $l$ was altered. For example, if $k=2$ in a $BN/BBN$, $k=2$ for $BH$ with $l=2$. While we do not believe the resulting $BH$ structures to be representative of real world hypergraph data sets - each element of $\mathcal{E}$ has in-degree = $l$, out-degree = 1 - we wanted to decouple the effects of changing $l$ from those of changing $k$ (i.e. we sought the most minimal changes between $BBN$ and $BH$).

\begin{figure}[ht]
\centering
\includegraphics[scale=0.225]{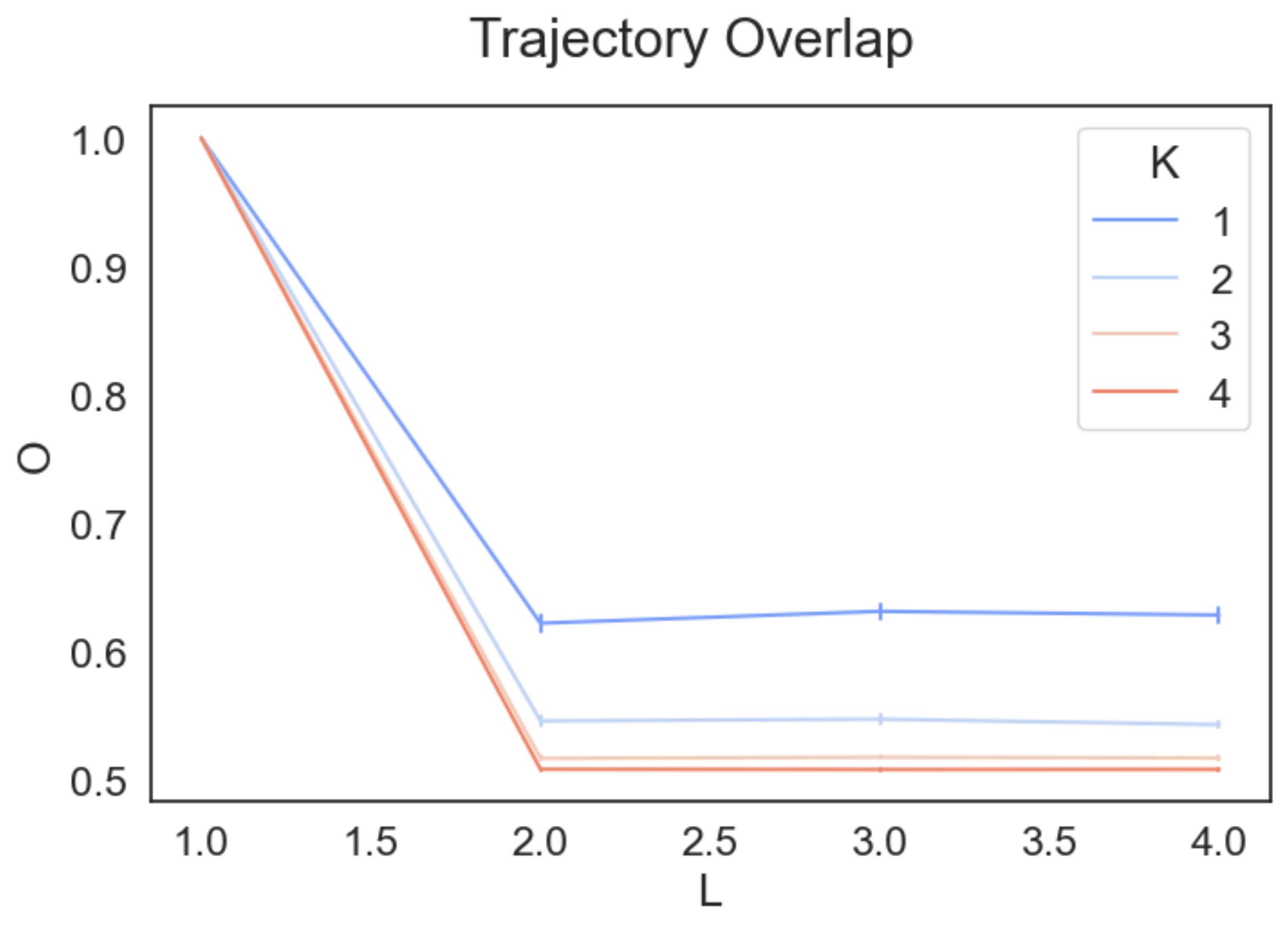}
\caption{State space trajectory overlap. Ensembles of 200 $BNs$ were constructed, and their derived $BBNs$ and $BHs$ simulated for 400 time steps, for each pair of $l$ and $k$. Mean trajectory overlaps for all vertices are reported for each condition.}
\label{fig:trajectory}
\end{figure}

Figure \ref{fig:trajectory} shows trajectory overlap results for $k$ and $l$ ranging from 1 to 4. For each condition, the ensemble mean is reported on networks that were simulated for 400 time steps. As expected, $O = 1$ when $l = 1$ for all values of $k$, and a decrease in $O$ is observed when $l$ is increased to 2. The magnitude of this decrease is amplified as $k$ moves from 1 to 4, indicating that increasing $l$ for higher values of $k$ results in increasingly non-overlapping state space trajectories. In the case of $k=4, l=2$, trajectories appear nearly independent. Interestingly, increasing $l$ beyond 2 does not seem to have drastic impacts on $O$ for all values of $k$. This indicates that increasing $l$ from 1 to 2 is all that is required to dramatically change state space trajectories between $BBNs$ and $BHs$. These results prompted us to examine how varying $l$ might change other characteristics of network dynamics such as attractor period lengths. 

\section{Results} \label{results}
Examining the attractor lengths of $BNs$ is a crucial method in understanding their state space dynamics \cite{borriello_basis_2021}. The attractor of a $BN$ is a periodic set of states that the network will enter - and will not exit unless it is perturbed - after a certain period of time. Attrators with period 1 are called fixed point or steady state attractors, an example of which can be seen in Figure \ref{fig:transition dynamics}A. Cyclic attractors have period $> 1$ and can vary from relatively simple behavior, like switching between two states, to the longer and more complex patterns observed in Figure \ref{fig:transition dynamics}C. Chaotic attractors are associated with much longer periods than those observed in Figure \ref{fig:transition dynamics} and an apparent lack of coherent trajectory in state space. Both Kauffman's original work and numerical analysis by Derrida and Pomeau \cite{derrida_random_1986} have revealed that $k=2$ represents a critical boundary between the ordered dynamics of fixed point and simple cyclic attractors on one hand, and chaotic dynamics on the other. It is therefore valuable to examine how attractor period lengths change when comparing $BHs$ to $BBNs$. 

\begin{figure}[ht]
\centering
\includegraphics[scale=0.25]{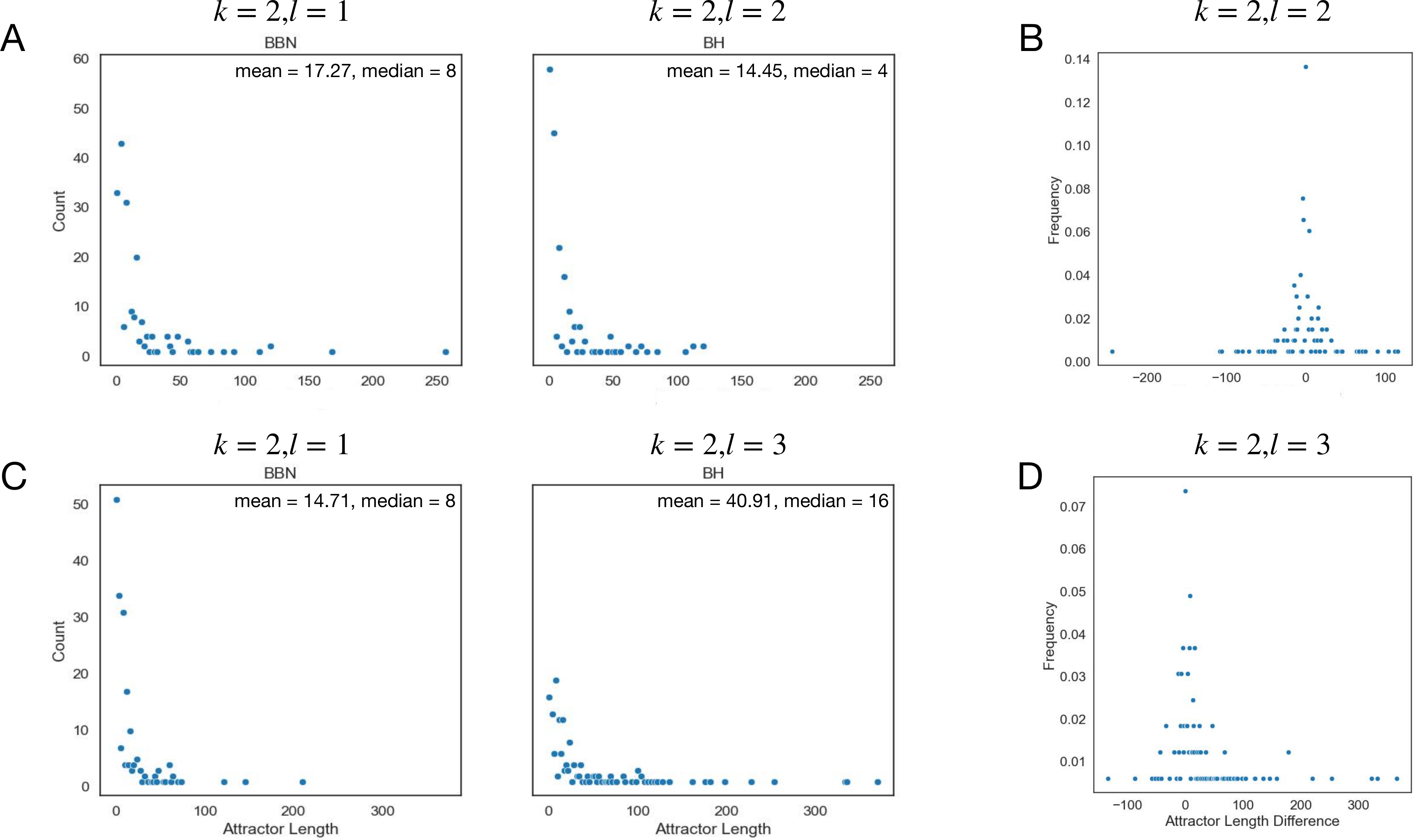}
\caption{Attractor period lengths of $BBNs$ and $BHs$ at $k=2$. (A) Attractor period length distribution for ensembles of randomly constructed $BBNs$ (left) and their derived $BHs$ for $l=2$ (right). (B) Distribution of the differences in attractor lengths between $BBNs$ and their derived $BH$ for $l=2$. Each point represents an individual experiment in which a $BH$ was derived from a $BBN$ and their attractor period lengths measured. Positive values along the x-axis indicate that $BH$ attractor length was greater than $BBN$ attractor length, and vice versa for negative values. (C) Attractor period length distributions for the case of comparing $BBNs$ (left) to $BHs$ with $l=3$ (right). For each $BH$ condition tested a unique ensemble of $BBNs$ was used, so the $BBN$ distribution shown here differs slightly from that in A. (D) Distribution of the differences in attractor lengths between $BBNs$ and their derived $BH$ for $k=2, l=3$.} 
\label{fig:attractors}
\end{figure}

Because $k=2$ is the critical boundary between ordered and chaotic dynamics, we decided to determine how attractor period lengths changed for $k=2$ when $l$ is increased beyond 1. Figure \ref{fig:attractors}A shows results obtained by constructing 200 random $BNs$ for $n=50$ and their corresponding $BBNs$ and $BHs$, with each network simulated for 400 times steps. Increasing $l$ from 1 (Figure \ref{fig:attractors}A, left) to 2 (Figure \ref{fig:attractors}A, right) yielded a modest decrease in mean period length from 17.27 to 14.15, although this result did not reach statistical significance. Median period length was also decreased from 8 to 4.  Interestingly, while summary statistics showed a slight decrease in period length, the change in period length was shown to vary drastically from experiment to experiment, with individual experiments yielding results from period length reduction by ~250 to period length increase by ~110 (Figure \ref{fig:attractors}B).

Increasing $l$ to 3 yielded quite different results, showing a marked increase in period length (Figure \ref{fig:attractors}C). In this case, mean attractor length was increased from 14.71 to 40.91 ($p=8.68 \times 10^{-8}$) as $l$ was increased from 1 (Figure \ref{fig:attractors}C, left) to 3 (Figure \ref{fig:attractors}C, right). Median period length was also doubled from 8 to 16. While increasing $l$ to 3 could still result in decreases in period length for a single experiment, the distribution was skewed towards increases (Figure \ref{fig:attractors}D). These results indicate that at the critical point between ordered and chaotic dynamics for $BNs$, attractor period lengths can be influenced by increases in network connectivity associated with moving from a $BBN$ to a $BH$. Furthermore, these changes depend on the extent to which connectivity increases, with $l=2$ yielding a slight reduction in mean and median period length and $l=3$ leading to substantial increases in both measures. The changes in attractor period length between $BBNs$ and $BHs$ at the critical boundary prompted us to probe for further differences in dynamical properties of these networks.

The complexity of a $BN$ is a measure of how well the network balances stability (ordered dynamics) and flexibility (chaotic dynamics) \cite{pineda_novel_2019}. In living systems, it is thought that this balance is important for the retention of useful information while also allowing an organism to explore new alternatives. $BN$ complexity is calculated by first computing the Shannon information entropy $E_i$ for each vertex $v_i, 1 \leq i \leq n$

\begin{equation}
    E_i = -(p_0 \log_2(p_0) + p_1 \log_2(p_1))
\end{equation}
where $p_0$ and $p_1$ are the proportions of time spent in state 0 and state 1, respectively, over the state space trajectory of the network. Letting $\bar{E} = \frac{1}{n} \sum^n_{i=1} E_i$, the complexity of the network is

\begin{equation}
    C = 4 \times \bar{E} \times (1 - \bar{E})
\end{equation}

\begin{figure}[ht]
\centering
\includegraphics[scale=0.25]{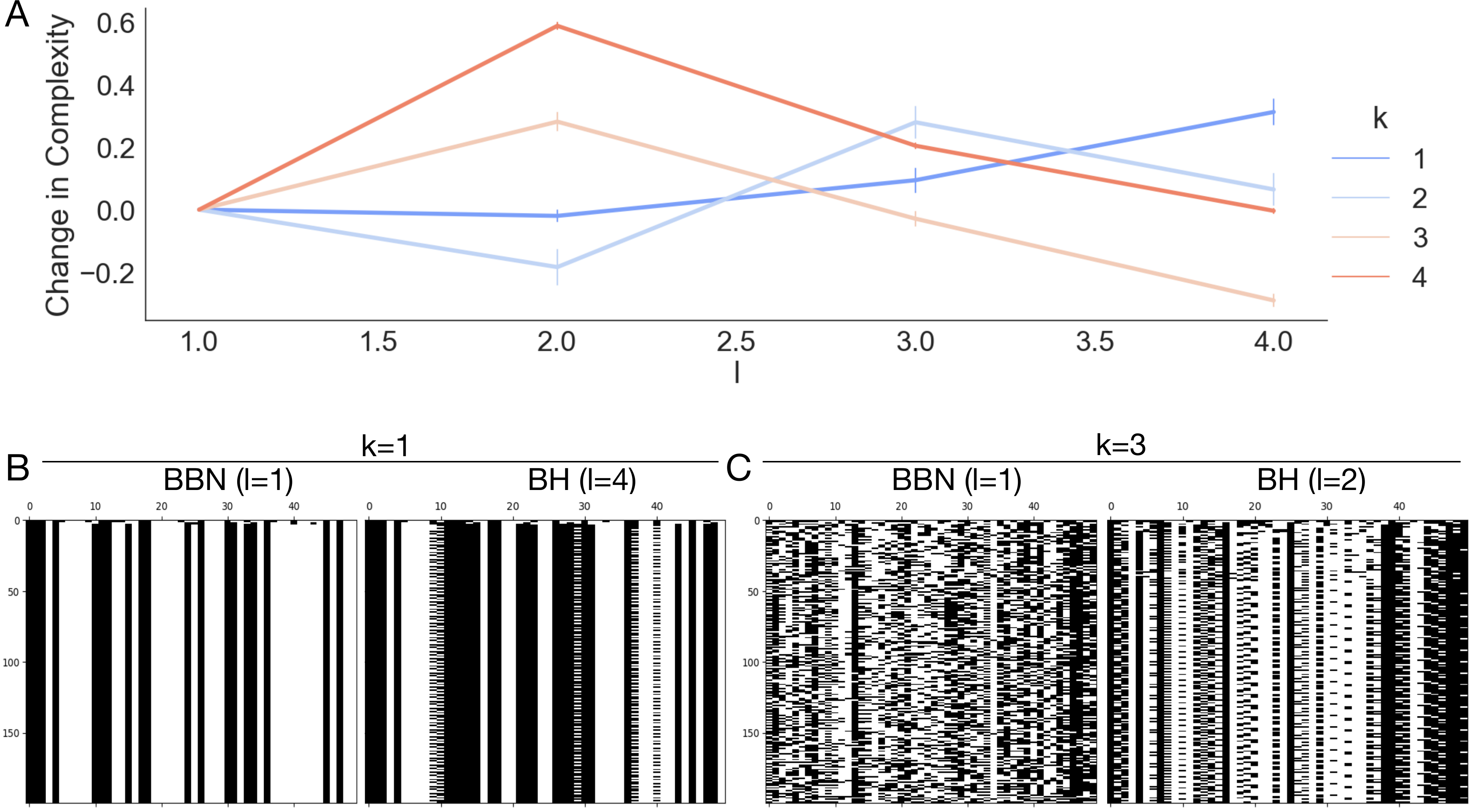}
\caption{Network complexity of $BBNs$ and $BHs$. (A) Mean difference in network complexity values between randomly constructed $BBNs$ and the $BHs$ derived from them for $k$ and $l$ values ranging from 1-4. Positive values indicate increased network complexity in BHs compared to BBNs. (B) Example of a change in state space trajectory brought about by increasing $l$ from 1 to 4 when $k=1$. (C) Example of a change in state space trajectory brought about by increasing $l$ from 1 to 2 when $k=3$.}
\label{fig:complexity}
\end{figure}

In order to determine if increasing $l$ leads to changes in network complexity, we simulated $BBNs$ and $BHs$ derived from $BNs$ of $n=100$. Values of $k$ and $l$ ranged from 1 to 4 and means of 200 independent simulations, each with 800 time steps, per condition are reported. As indicated by Figure \ref{fig:complexity}A, increasing $l$ tends to lead to increases in $BH$ complexity compared to $BBN$ controls, with the opposite only occurring in the $k=2, l=2$ and $k=3, l=4$ cases. A hypothesis for the increase in complexity as $l$ increases for $k=1$ is that the addition of new edges in $BHs$ pushes state space dynamics from the ordered phase towards criticality. In agreement with this hypothesis, average complexity for $k=1, l=4$ begins to approach the complexity values observed for $BBNs$ when $k=2$.

Recall that for a $BN$ consisting of $n$ vertices, the $BBN$ and $BH$ derived from it will contain $n$ $\mathcal{V}$ vertices and $kn$ $\mathcal{E}$ vertices. This implies that as $k$ increases, measures of state space dynamics like network complexity will largely be driven by the $\mathcal{E}$ vertices. Even in the case where $k=1$, the $BBN$ and $BH$ derived from a $BN$ will contain $50 \%$ each of $\mathcal{V}$ vertices and $\mathcal{E}$ vertices. To validate our hypothesis above and ensure that our findings are not solely a result of changes in the dynamics of $\mathcal{E}$ vertices, we examined state space trajectories of $\mathcal{V}$ vertices on time steps where this subset was updated. We observed that $\mathcal{V}$ vertices in $BBNs$ tend towards fixed point attractors while their corresponding $BHs$ tend towards cyclic attractors as $l$ increases, an example of which is shown in Figure \ref{fig:complexity}B for the case where $l=4$.

Mean $BH$ complexity is reduced compared to $BBN$ controls for $k=2, l=2$. Interestingly, while $BH$ complexity fell compared to $BBN$ controls for $k=2, l=2$, an increase is observed for $k=3, l=2$ (Figure \ref{fig:complexity}A). One potential reason is that increasing $l$ to 2 for $k=3$ could have stabilizing effects on state space dynamics. In support of this hypothesis, several examples of state space trajectories were found in which what appear to be chaotic attractors in $BBNs$ were driven to cyclic attractors when $l$ was increased (Figure \ref{fig:complexity}C). Additionally, we observed that mean $BH$ attractor length was 49, compared to 54.17 for $BBN$ controls, confirming that the increase in complexity for $k=3, l=2$ is accompanied by a decrease in mean attractor period length. This provides further evidence that increasing $l$ to 2 when $k=3$ has stabilizing effects on state space trajectories. When examining the results for $k=4$, we observe complexity dynamics similar to what was observed for $k=3$ as $l$ is increased: $BH$ complexity increased drastically from $l=1$ to $l=2$ before falling again as $l$ increased further. While we believe the reasons for this are similar to the $k=3$ case, attractor lengths were often too long to be measured in the time frame each network was simulated in, leaving this as an open question for future work.

Thus far we have analyzed simulations in which network dynamics were allowed to proceed unperturbed after initial states were set. While valuable, the insights gained from these experiments neglect an oft-examined feature of Boolean networks - how they respond to random perturbations \cite{baudin_controlling_2019}. To elucidate how $BHs$ respond to perturbations, we employed a measure of fragility previously applied to $BNs$. Briefly, we measured the mean change in network fragility between a $BH$ and the $BBN$ that extends it over 200 independent experiments as described in \cite{pineda_novel_2019}. During each experiment, we randomly chose 20 vertices at each of 400 time steps and changed their state variable. The value of 20 was chosen because it was shown to induce differences in network fragility for values of $k$ ranging from 1 to 5 in \cite{pineda_novel_2019}.

\begin{figure}[ht]
\centering
\includegraphics[scale=0.24]{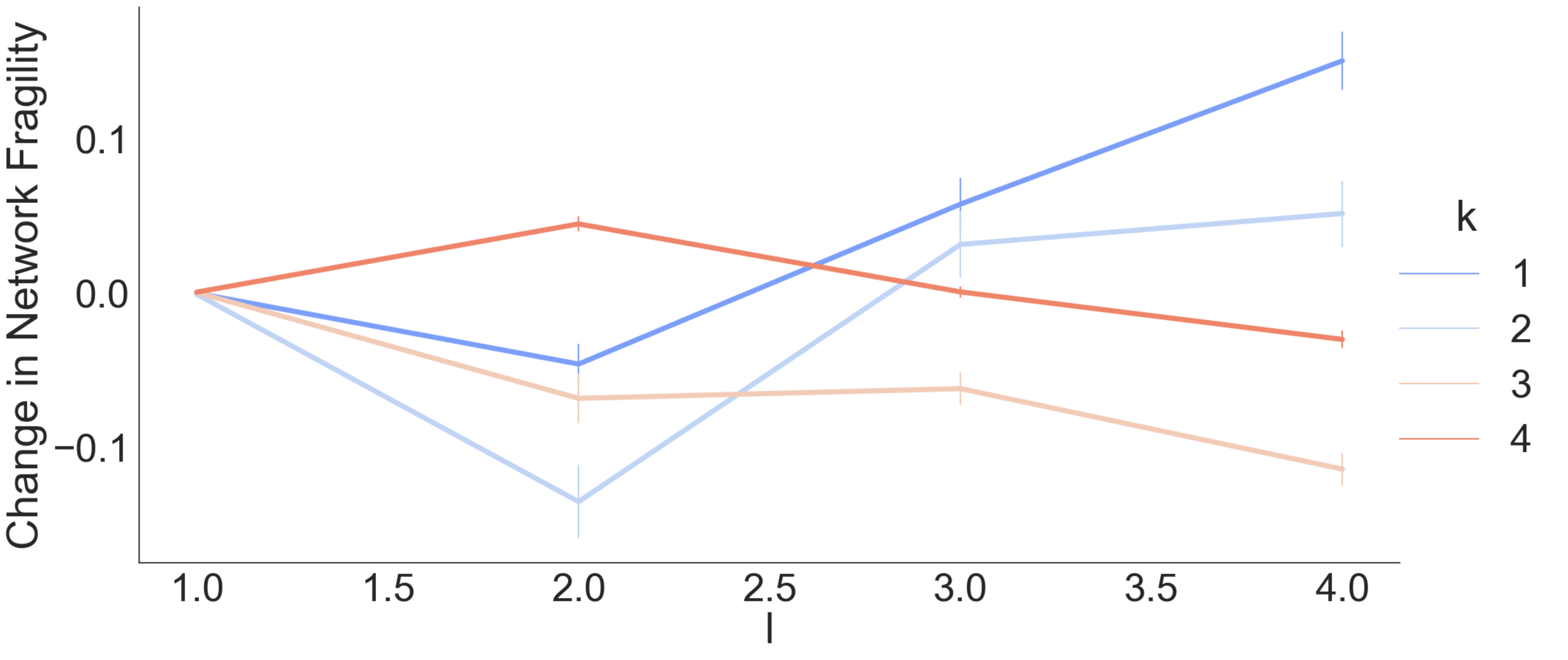}
\caption{Change in network fragility between $BHs$ and $BBNs$. $BBNs$ and $BHs$ derived from $BNs$ of $n=50$ were subject to perturbed simulations as described in the text. The mean difference in network fragility for ensembles of 200 independent experiments is reported for each condition. Negative values indicate that $BHs$ are more antifragile than their $BBN$ controls, and vice versa for positive values.}
\label{fig:fragility}
\end{figure}

Our results indicate changes in network fragility between $BHs$ and their $BBN$ controls across all $k$ from 1 to 4, the magnitude of which varying depending on the value of $l$ (Figure \ref{fig:fragility}). For $k=1$, increasing $l$ lead to steady increases in mean fragility for $l > 2$. For $k=2, l=2$, $BHs$ were more antifragile than their corresponding $BBNs$, but this decrease vanished as $l$ increased beyond 2. For both $k=3$ and $k=4$, increases in $l$ tended to lead to increased antifragility in $BHs$ compared to $BBNs$.

\section{Discussion} \label{discussion}
Our results indicate that random Boolean hypernetworks exhibit differences in state space behavior as measured by attractor period length, network complexity and network fragility when compared to random Boolean networks. In some cases, the higher order interactions present in $BHs$ stabilized state space trajectory compared to $BBN$ controls (Figure \ref{fig:complexity}C), while in other cases $BH$ dynamics appeared more chaotic (Figure \ref{fig:attractors}C). These results indicate that the effect of higher order interactions on Boolean network dynamics are not generally stability-inducing or instability-inducing, but instead depend on the specific values of $k$ and $l$ chosen. Furthermore, the results presented in this study represent tendencies captured in ensemble means and do not indicate that higher order interactions will always be stability-inducing or instability-inducing for randomly assembled network with specific parameters $k$ and $l$.

The context-dependent effects on state space stability induced by higher order interactions in our model seem to align well with results obtained by Gallo, et al. in a hypergraph system of coupled Rossler oscillators \cite{gallo_synchronization_2022}. These results indicate that stabilizing or destabilizing effects of higher order interactions depend on the coupling strengths of higher order interactions. Higher order interactions with weak coupling strengths induce instability and vice versa, in line with the results obtained in our study. By this, we mean that increasing $l$ for low $k$ tends to lead from stable to critical or chaotic dynamics, while increasing $l$ for higher values of $k$ tends to shift networks from chaotic to critical dynamics. Building on the work of Gallo, Li, et al. built a model of social contagion spread on hypernetworks \cite{li_enhancing_2024}. While similar to our study in terms of describing state dynamics on hypergrpahs, Li, et al. used a mean field approach to determine contagion spread.

It should be noted that all hypernetworks in our study were composed of hyperedges with a single vertex in their heads. For cases where $l \ge 2$, such hyperedges correspond to B-edges described in \cite{gallo_directed_1993}, meaning that all networks where $l \ge 2$ were B-networks. While using B-networks allowed us to interrogate the effects of increasing $l$ for specific values of $k$, real world hypergraph data is unlikely to be structured as B-networks. Therefore, we believe that future work on Boolean hypernetworks should focus in part on modeling real world systems in which hyperedges contain multiplicity of inputs and outputs. Generalizing in this way will require us to carefully consider how we treat self-loops. In this work, we depart from \cite{gallo_directed_1993} by loosening the requirement that each directed hyperedge contain disjoint subsets of vertices to allow for hypergraph self-loops. While our desire to include self-loops in our model stems from the importance of auto-regulation in biological systems \cite{roy_autoregulation_2020}, care must be taken to ensure that we generalize these structures beyond the B-networks described in this work.

Further work on Boolean hypernetworks should also address the mathematical foundations of these systems, especially with respect to the relationships between Boolean networks, bipartite Boolean networks, and Boolean hypernetworks. For instance, the sequence of transformations outlined in Figure \ref{fig:network relationships} to transform a $BN$ to a $BBN$ require both a transformation of network architecture (from a directed graph to a directed bipartite graph) and a mapping of state update rules that faithfully recapitulates the behavior of the $BN$ (i.e. $\mathcal{E}$ edges must take the identity rule). Transforming a $BBN$ to a $BH$ also requires both types of transformations: additional edges that terminate in $\mathcal{E}$ vertices must be added, and $\mathcal{E}$ vertices can no longer be assigned the identity rule for $l > 2$. Such mappings in which transformations of both structure and function are required lend themselves to category theoretical analysis \cite{Spivak2014CategoryTF}, which should be utilized in future work. Category theory has already been applied to graphs \cite{guerra_using_2022, haruna_application_nodate}, and the extension of this application to the relationship between $BNs$ and $BHs$ could yield valuable insights.

Finally, the networks generated and analyzed in this study were all derived from $BNs$ that contained either 50 or 100 vertices. While we observed strong similarities in results between networks derived from 50 vertex $BNs$ and 100 vertex $BNs$ across all experiments, future work should be aimed at examining the dynamical properties of $BHs$ with a wider variety of network topologies. Along with examining larger and smaller networks, with respect to vertex number, for fixed $k$ and $l$, randomly generated $BHs$ with scale free toplogies should also be examined. Methods for generating hypernetworks with scale free topologies already exist in the literature and could be used to further this aim \cite{aksoy_hypernetwork_2020}.

\bibliographystyle{plain}
\bibliography{refs}

\end{document}